\documentclass[12pt, reqno]{amsart}
\usepackage{amssymb,amsthm,amsfonts,amsmath, amscd}

\usepackage[dvips]{graphicx}

\usepackage{hyperref}
\usepackage{mathrsfs}

\newcommand{\Ga}{\Gamma}

\newcommand{\be}{\beta}

\newcommand{\la}{\lambda}

\newcommand{\tth}{\theta}

\newcommand{\La}{\Lambda}
\newcommand{\epsi}{\varepsilon}
\newcommand{\De}{\Delta}
\newcommand{\si}{\sigma}
\newcommand{\Om}{\Omega}

\newcommand{\const}{\operatorname{const}}

\newcommand{\R}{\mathbb R}
\newcommand{\C}{\mathbb C}
\newcommand{\Z}{\mathbb Z}
\newcommand{\N}{\mathbb N}

\newcommand{\Y}{\mathbb Y}

\newcommand{\Part}{\mathscr P}

\newcommand{\PP}{\mathcal P}
\newcommand{\QQ}{\mathcal Q}

\newcommand{\Sym}{\mathfrak S}
\newcommand{\HH}{{\mathfrak H}}
\newcommand{\UU}{{\mathfrak U}}

\newcommand{\wt}{\widetilde}
\newcommand{\wh}{\widehat}
\newcommand{\n}{{(n)}}
\newcommand{\zz}{{z,z'}}
\newcommand{\zzxi}{{z,z',\xi}}

\newcommand{\re}{\operatorname{Re}}

\newtheorem{theorem}{Theorem}[section]

\newtheorem{lemma}[theorem]{Lemma}

\theoremstyle{definition}

\newtheorem{remark}[theorem]{Remark}

\numberwithin{equation}{section}

\begin{document}

\title[Random permutations and related topics]
{Random permutations and related topics}

\author{Grigori Olshanski}
\address{Institute for Information Transmission Problems, Bolshoy Karetny 19, Moscow 127994, Russia;
\newline\indent Independent University of Moscow, Russia} \email{olsh2007@gmail.com}

\date{}
\begin{abstract}
We present an overview of selected topics in random permutations and random
partitions highlighting analogies with random matrix theory.
\end{abstract}

\maketitle

\tableofcontents

\section{Introduction}\label{1}

An ensemble of random permutations is determined by a probability distribution
on $S_n$, the set of permutations of $[n]:=\{1,2,\dots,n\}$. Even the simplest
instance of the uniform probability distribution is already very interesting
and leads to a deep theory. The symmetric group $S_n$ is in many ways linked to
classical matrix groups, and  ensembles of random permutations should be
treated on equal footing with random matrix ensembles, such as the ensembles of
classical compact groups and symmetric spaces of compact type with the
normalized invariant measure. The role of matrix eigenvalues is then played by
partitions of $n$ that parameterize the conjugacy classes in $S_n$. The
parallelism with random matrices becomes especially striking in applications to
constructing representations of ``big groups'' --- inductive limits of
symmetric or classical compact groups.

The theses stated above are developed in Section \ref{2}. The two main themes
of this section are the space $\Sym$ of virtual permutations ($\Sym$ is a
counterpart of the space of Hermitian matrices of infinite size) and the
Poisson--Dirichlet distributions (a remarkable family of infinite-dimensional
probability distributions). We focus on a special family of probability
distributions on $S_n$ with nice properties, the so--called Ewens measures
(they contain the uniform distributions as a particular case). It turns out
that the large-$n$ limits of the Ewens measures can be interpreted as
probability measures on $\Sym$. On the other hand, the Ewens measures give rise
to ensembles of random partitions, from which one gets, in a limit transition,
the Poisson--Dirichlet distributions.

A remarkable discovery of Frobenius, the founder of representation theory, was
that partitions of $n$ not only parameterize conjugacy classes in $S_n$ but
also serve as natural coordinates in the dual space $\wh S_n$
--- the set of equivalence classes of irreducible representations of $S_n$.
This fact forms the basis of a ``dual'' theory of random partitions, which
turns out to have many intersections with random matrix theory. This is the
subject of Sections \ref{3}--\ref{4}. Here we survey results related to the
Plancherel measure on $\wh S_n$ and its consecutive generalizations: the
z-measures and the Schur measures.

Thus, the two-faced nature of partitions gives rise to two kinds of
probabilistic models. At first glance, they seem to be weakly related, but
under a more general approach one sees a bridge between them. The idea is that
the probability measures in the ``dual picture'' can be further generalized by
introducing an additional parameter, which is an exact counterpart of the $\be$
parameter in random matrix theory. This parameter interpolates between the
group level and the dual space level, in the sense that in the limit as
$\be\to0$,  the ``beta z-measures'' degenerate to the measures on partitions
derived from the Ewens measures, see Section \ref{4-4} below.

\section[The Ewens measures]{The Ewens measures, virtual permutations, and the Poisson--Dirichlet
distributions}\label{2}

\subsection{The Ewens measures}\label{2-1}

Permutations $s\in S_n$ can be represented as $n\times n$ unitary matrices
$[s_{ij}]$, where $s_{ij}$ equals 1 if $s(j)=i$, and 0 otherwise. This makes it
possible to view random permutations as a very special case of random unitary
matrices.

Given a probability distribution on each set $S_n$, $n=1,2,\dots$, one may
speak of a sequence of ensembles of {\it random permutations\/} and study their
asymptotic properties as $n\to\infty$. The simplest yet fundamental example of
a probability distribution on $S_n$ is the {\it uniform distribution\/} $P^\n$,
which gives equal weights $1/n!$ to all permutations $s\in S_n$; this is also
the normalized Haar measure on the symmetric group. However, a more complete
picture is achieved by considering a one-parameter family of distributions,
$\{P^\n_\tth\}_{\tth>0}$, forming a deformation of $P^\n$:
\begin{equation}\label{ewensmeasure}
P^\n_\tth(s)=\frac{\tth^{\ell(s)}}{\tth(\tth+1)\dots(\tth+n-1)}, \qquad s\in
S_n,
\end{equation}
where $\ell(s)$ denotes the number of cycles in $s$ (the uniform distribution
corresponds to $\tth=1$).

For reasons explained below we call $P^\n_\tth$ the {\it Ewens measure\/} on
the symmetric group $S_n$ with parameter $\tth$. Obviously, the Ewens measures
are invariant under the action of $S_n$ on itself by conjugations.

We propose to think of $(S_n,P^\n)$ as of an analogue of CUE$_N$, Dyson's
circular unitary ensemble \cite{For10a} formed by the unitary group $U(N)$
endowed with the normalized Haar measure $P^{U(N)}$. More generally, the Ewens
family $\{P^\n_\tth\}$ should be viewed as a counterpart of a family of
probability distributions on the unitary group $U(N)$ forming a deformation of
the Haar measure:
\begin{equation}
P^{U(N)}_t(dU)=\const\, |\det((1+U)^t)|^2 P^{U(N)}(dU), \qquad U\in U(N),
\end{equation}
where $t$ is a complex parameter, $\re t>-\frac12$ (the Haar measure
corresponds to $t=0$). Using the Cayley transform one can identify the manifold
$U(N)$, within a negligible subset, with the flat space $H(N)$ of $N\times N$
Hermitian matrices; then $P^{U(N)}_t$ turns into the so-called {\it
Hua--Pickrell measure\/} on $H(N)$:
\begin{equation}
P^{H(N)}_t(dX)=\const\,|\det(1+iX)^{-t-N}|^2 dX,
\end{equation}
where ``$dX$'' in the right-hand side denotes Lebesgue measure. For more
detail, see \cite{Bor01b}, \cite{Ner02}, \cite{Ols03a}. A similarity between
$P^\n_\tth$ and $P^{U(N)}_t$ (or $P^{H(N)}_t$) is exploited in \cite{Bou07}.

A fundamental property of the Ewens measures is their consistency with respect
to some natural projections $S_n\to S_{n-1}$ that we are going to describe:

Given a permutation $s\in S_n$, the {\it derived permutation\/} $s'\in S_{n-1}$
sends each $i=1,\dots,n-1$ either to $s(i)$ or to $s(n)$, depending on whether
$s(i)\ne n$ or $s(i)=n$. In other words, $s'$ is obtained by removing $n$ from
the cycle of $s$ that contains $n$. For instance, if $s=(153)(24)$ (meaning
that one cycle in $s$ is $1\to5\to3\to1$ and the other is $2\to4\to2$), then
$s'=(13)(24)$. The map $S_n\to S_{n-1}$ defined in this way is called the {\it
canonical projection\/} \cite{Ker04} and denoted as $p_{n-1,n}$. For $n\ge5$,
it can be characterized as the only map $S_n\to S_{n-1}$ commuting with the
two-sided action of the subgroup $S_{n-1}$.

The next assertion \cite{Ker04} is readily verified:

\begin{lemma}
For any $n=2,3,\dots$, the push-forward of $P^\n_\tth$ under $p_{n-1,n}$
coincides with $P^{(n-1)}_\tth$.
\end{lemma}

The Hua--Pickrell measures enjoy a similar consistency property with respect to
natural projections $p^H_{N-1,N}\colon H(N)\to H(N-1)$ (removal of the $N$th
row and column from an $N\times N$ matrix) \cite{Bor01b}, \cite{Ols03a}. This
fact is hidden in the old book by Hua Lokeng \cite{Hua58}, in his computation
of the matrix integral
\begin{equation}\label{hua}
\int_{H(N)}\det(1+X^2)^{-t}dX
\end{equation}
by induction on $N$. (Note that \cite{Hua58} contains a lot of masterly
computations of matrix integrals.) Much later, the consistency property was
rediscovered and applied to constructing measures on infinite-dimensional
spaces by Shimomura \cite{Shi75}, Pickrell \cite{Pic87}, and Neretin
\cite{Ner02}. Note that analogues of the projections $p^H_{N-1,N}$ can be
defined for other matrix spaces including the three series of compact classical
groups and, more generally, the ten series of classical compact symmetric
spaces \cite{Ner02}.

\subsection{Virtual permutations, central measures, and Kingman's
theorem}\label{2-2}

As we will see, the consistency property of the Ewens measures $P^\n_\tth$
makes it possible to build an ``$n=\infty$'' version of these measures. In the
particular case $\tth=1$, this leads to a concept of  ``uniformly distributed
infinite permutations''.

Let $\Sym=\varprojlim S_n$ be the projective limit of the sets $S_n$ taken with
respect to the canonical projections. By the very definition, each element
$\si\in\Sym$ is a sequence $\{\si_n\in S_n\}_{n=1,2,\dots}$ such that
$\si_{n-1}=p_{n-1,n}(\si_n)$ for any $n=2,3,\dots$\,. We call $\Sym$ the {\it
space of virtual permutations\/} \cite{Ker04}. It is a compact topological
space with respect to the projective limit topology.

By classical Kolmogorov's theorem, any family $\{\mathcal P^\n\}$ of
probability measures on the groups $S_n$, consistent with the canonical
projections, gives rise to a probability measure $\mathcal
P=\varprojlim\mathcal P^\n$ on the space $\Sym$. Taking $\mathcal P^\n=
P^\n_\tth$, $\tth>0$, we get some measures $P_\tth$ on $\Sym$ with nice
properties; we still call them the Ewens measures.

A parallel construction exists in the context of matrix spaces. In particular,
by making use of the projections $p^H_{N-1,N}$, one can define a projective
limit space $\HH:=\varprojlim H(N)$, which is simply the space of all Hermitian
matrices of infinite size. This space carries the measures
$P^\HH_t:=\varprojlim P^{H(N)}_t$, $\re t>-\frac12$. Using the Cayley
transform, the measures $P^\HH_t$ can be transformed to some measures $P^\UU_t$
on a projective limit space $\UU:=\varprojlim U(N)$.

As is argued in \cite{Bor01b}, the $t=0$ case of the probability space
$(\HH,P^\HH_t)\simeq(\UU, P^\UU_t)$ may be viewed as an ``$N=\infty$'' version
of CUE$_N$. Likewise, we regard $(\Sym, P_1)$ as an ``$n=\infty$'' version of
the finite uniform measure space $(S_n,P^\n_1)$.

A crucial property of CUE$_N$ is its invariance under the action of the unitary
group $U(N)$ on itself by conjugation. This property is shared by the deformed
ensemble with parameter $t$ as well. In the infinite-dimensional case, the
analogous property is $U(\infty)$-invariance of the measures $P^\HH_t$; here by
$U(\infty)$ we mean the inductive limit group $\varinjlim U(N)$, which acts on
the space $\mathfrak H$ by conjugations.

Likewise, define the {\it infinite symmetric group\/} $S_\infty$ as the
inductive limit $\varinjlim S_n$, or simply as union of the finite symmetric
groups $S_n$, where each group $S_n$ is identified with the subgroup of
$S_{n+1}$ fixing the point $n+1$ of the set $[n+1]$. The group $S_\infty$ is
countable and can be realized as the group of all permutations of the set
$\N=\{1,2,\dots\}$ moving only finitely many points. The space $\Sym$ contains
$S_\infty$ as a dense subset, and the action of the group $S_\infty$ on itself
by conjugation extends by continuity to the space $\Sym$; we will still call
the latter action the {\it conjugation action\/} of $S_\infty$.

It is convenient to give a name to measures that are invariant under the
conjugation action of $S_n$ or $S_\infty$; let us call such measures {\it
central\/} ones. Now, a simple lemma says:

\begin{lemma}\label{central}
Let, as above, $\mathcal P=\varprojlim\mathcal P^\n$ be a projective limit
probability measure on $\Sym$. If $\mathcal P^\n$ is central for each $n$, then
so is $\mathcal P$.
\end{lemma}

As a consequence we get that the measures $P_\tth$ are central.

In a variety of random matrix problems, the invariance property under an
appropriate group action makes it possible to pass from matrices to their
eigenvalues or singular values. For random permutations from $S_n$ directed by
a central measure, a natural substitute of eigenvalues is another invariant
--- partitions of $n$ parameterizing the cycle structure of permutations.

In combinatorics, by a {\it partition\/} one means a sequence
$\rho=(\rho_1,\rho_2,\dots)$ of weakly decreasing nonnegative integers with
infinitely many terms and finite sum $|\rho|:=\sum \rho_i$. Of course, the
number of nonzero terms $\rho_i$ in $\rho$ is always finite; it is denoted by
$\ell(\rho)$. The finite set of all partitions $\rho$ with $|\rho|=n$ will be
denoted as $\Part(n)$. To a permutation $s\in S_n$ we assign a partition
$\rho\in\Part(n)$ (in words, a partition of $n$) comprised by the cycle-lengths
of $s$ written in weakly decreasing order. Obviously, $\rho$ is a full
invariant of the conjugacy class of $s$. The projection $s\mapsto\rho$ takes
any probability measure on $S_n$ to a probability measure on $\Part(n)$. Now,
the point is that this establishes a one-to-one correspondence $\mathcal
P^\n\leftrightarrow\Pi^\n$ between arbitrary central probability measures
$\mathcal P^\n$ on $S_n$ and arbitrary probability measures $\Pi^\n$ on
$\Part(n)$. In this sense, random permutations $s\in S_n$ (directed by a
central measure) may be replaced by {\it random partitions\/}
$\rho\in\Part(n)$.

The link between $\mathcal P^\n$ and $\Pi^\n$ is simple: Given
$\rho\in\Part(n)$, let $C(\rho)\subset S_n$ denote the corresponding conjugacy
class in $S_n$. Then $\Pi^\n(\rho)=|C(\rho)|\mathcal P^\n(s)$ for any $s\in
C(\rho)$. Further, there is an explicit expression for $|C(\rho|$: it is equal
to $n!/z_\rho$, where $z_\rho=\prod k^{m_k}m_k!$ and $m_k$ stands for the
multiplicity of $k=1,2,\dots$ in $\rho$.

In this notation, the measure $\Pi^\n_\tth$ corresponding to the Ewens measure
$\mathcal P^\n=P^\n_\tth$ is given by the expression
\begin{equation}\label{ESF}
\Pi^\n_\tth (\rho)=\frac{\tth^{\ell(\rho)}}{\tth(\tth+1)\dots(\tth+n-1)}
\,\prod_k\frac1{k^{m_k}m_k!}\,, \qquad \rho\in\Part(n),
\end{equation}
widely known under the name of {\it Ewens sampling formula\/} \cite{Ewe98}.
This justifies the name given to the measure \eqref{ewensmeasure}.

The following result provides a highly nontrivial ``$n=\infty$'' version of the
evident correspondence $\mathcal P^\n\leftrightarrow\Pi^\n$:

\begin{theorem}\label{kingman}
There exists a natural one-to-one correspondence $\mathcal P\leftrightarrow\Pi$
between arbitrary central probability measures $\mathcal P$ on the space $\Sym$
of virtual permutations and arbitrary probability measures $\Pi$ on the space
\begin{equation}
\overline{\nabla}_\infty:=\{(x_1,x_2,\dots)\in[0,1]^\infty\colon x_1\ge
x_2\ge\dots, \quad \sum x_i\le1\}.
\end{equation}
\end{theorem}

In other words, each central measure $\mathcal P$ on $\Sym$ is uniquely
representable as a mixture of {\it indecomposable\/} (or {\it ergodic\/})
central measures, which in turn are parameterized by the points of
$\overline\nabla_\infty$. The measure $\Pi$ assigned to $\mathcal P$ is just
the mixing measure.

\begin{proof}[Idea of proof] The theorem is a reformulation of celebrated {\it
Kingman's theorem\/}, see \cite{Kin78b}. Kingman did not deal with virtual
permutations but worked with some sequences of random permutations that he
called {\it partition structures\/}. Represent $\mathcal P$ as a projective
limit measure, $\mathcal P=\varprojlim\mathcal P^\n$. By Lemma \ref{central},
all measures $\mathcal P^\n$ are central. Pass to the corresponding measures
$\Pi^\n$ on partitions. The consistency of the family $\{\mathcal P^\n\}$ with
the canonical projections $S_n\to S_{n-1}$ then translates as the consistency
of the family $\{\Pi^\n\}$ with some canonical Markov transition kernels
$\Part(n)\to\Part(n-1)$. In Kingman's language this just means that $\{\mathcal
P^\n\}$ is a partition structure. Kingman's theorem provides a kind of Poisson
integral representation of partition structures via probability measures on
$\overline{\nabla}$, which is equivalent to the claim of Theorem \ref{kingman}.
\end{proof}

\medskip

Other proofs of Kingman's theorem can be found in \cite{Ker03}, \cite{Ker98},
where this result is placed in the broader context of potential theory for
branching graphs. For our purpose it is worth emphasizing that the claim of
Theorem \ref{kingman} has a counterpart in the random matrix context ---
description of $U(\infty)$-invariant probability measures on $\mathfrak H$,
which in turn is equivalent to classical Schoenberg's theorem on totally
positive functions \cite{Sch51}, \cite{Pic91}, \cite{Ols96}.

{}From the proof of Kingman's theorem it is seen that the space
$\overline\nabla_\infty$ arises as a large-$n$ limit of finite sets $\Part(n)$,
and every measure $\Pi$ can be interpreted as a limit of the corresponding
measures $\Pi^\n$. In this picture, the ergodic measures
$\Pi_{(x_1,x_2,\dots)}$ that are parameterized by points
$(x_1,x_2,\dots)\in\overline\nabla_\infty$ arise as limits of uniform
distributions on conjugacy classes $C(\rho)$ in growing finite symmetric
groups. Thus, it is tempting to regard the $\Pi_{(x_1,x_2,\dots)}$'s as a
substitute of those uniform measures.

\subsection{Application to representation theory}\label{2-3}

Pickrell's pioneer work \cite{Pic87} demonstrated how some non-Gaussian
measures on infinite-dimensional matrix spaces can be employed in
representation theory. The idea of \cite{Pic87} was further developed in
\cite{Ols03a}. As shown there, the measures $P^\HH_t$ on $\HH$ (or,
equivalently, the measures $P^\UU_t$ on $\UU$) meet a lack of the Haar measure
in the infinite-dimensional situation and can be applied to the construction of
some ``generalized (bi)regular representations'' of the group $U(\infty)\times
U(\infty)$.

Pickrell's work was also a starting point for a parallel theory for the group
$S_\infty$, \cite{Ker93c}, \cite{Ker04}. The key point is that the Ewens
measures $P_\tth$ on $\Sym$ have good transformation properties with respect to
a natural action of the group $S_\infty\times S_\infty$ on $\Sym$ extending the
two-sided action of $S_\infty$ on itself. Namely, the measure $P_1$ is
$S_\infty\times S_\infty$-invariant and is the only probability measure on
$\Sym$ with such a property, so that it may be viewed as a substitute of the
uniform distribution on $S_n$. As for the Ewens measures $P_\tth$ with general
$\tth>0$, they turn out to be quasi-invariant with respect to the action of
$S_\infty\times S_\infty$. The quasi-invariance property forms the basis of the
construction of some ``generalized (bi)regular representations'' $T_z$ of the
group $S_\infty\times S_\infty$. Here $z$ is a parameter ranging over $\C$, and
the Hilbert space of $T_z$ is $L^2(\Sym,P_{|z|^2})$. We refer to \cite{Ker04}
and \cite{Ols03b} for details.

\subsection{Poisson--Dirichlet distributions}\label{2-4}

The probability measures on $\overline\nabla_\infty$ assigned by Theorem
\ref{kingman} to the Ewens measures $P_\tth$ are known under the name of {\it
Poisson--Dirichlet distributions\/}; denote them by $PD(\tth)$. Continuing our
juxtaposition of the Ewens measures and the Hua--Pickrell measures one may say
that the Poisson--Dirichlet distributions are counterparts of the determinantal
point processes directing the decomposition of the measures $P^\HH_t$ into
ergodic components (those processes involve, in a slightly disguised form, the
sine-kernel process, see \cite{Bor01b}). Although the Poisson--Dirichlet
distributions $PD(\tth)$ seem to be much simpler than the sine-kernel process,
they are still very interesting objects with a rich structure. Below we list a
few equivalent descriptions of the $PD(\tth)$'s:

\smallskip

(a) {\it Projection of a Poisson process\/}. Let $P(\tth)$ denote the
inhomogeneous Poisson process on the half-line $\R_{>0}:=\{\tau\in\R\mid
\tau>0\}$ with intensity $\tth \tau^{-1}e^{-\tau}$, and let $y=\{y_i\}$ be the
random point configuration on $\R_{>0}$ with law $P(\tth)$. Due to the fast
decay of the intensity at $\infty$, the configuration $y$ is almost surely
bounded from above, so that we may arrange the $y_i$'s in weakly decreasing
order: $y_1\ge y_2\ge\dots>0$. Furthermore, the sum $r:=\sum y_i$ is almost
surely finite. Finally, it turns out that $r$ and the normalized vector
$x:=(y_1/r,y_2/r,\dots)\in\nabla_\infty$ are independent from each other, the
random variable $r$ has the gamma distribution on $\R_{>0}$ with density
$(\Ga(\tth))^{-1}t^{\tth-1}\exp(-t)$, and the random vector $x$ is distributed
according to $PD(\tth)$.

This means, in particular, that $PD(\tth)$ arises as the push-forward of the
Poisson process $P(\tth)$ under the projection $y\mapsto x$.

\medskip

(b) {\it Limit of Dirichlet distributions\/} \cite{Kin75}. Let $D_n(\tth) $
denote the probability distribution on the $(n-1)$-dimensional simplex
$$
\De_n:=\{(x_1,\dots, x_n)\in\R^n\mid x_1,\dots,x_n\ge0, \quad \sum x_i=1\}
$$
with the density proportional to $\prod x_i^{n^{-1}\tth-1}$ (with respect to
Lebesgue measure on $\De_n$). Note that $D_n(\tth)$ enters the family of the
Dirichlet distributions.  Rearranging the coordinates $x_i$ in weakly
decreasing order and adding infinitely many 0's one gets a map
$\De_n\to\overline{\nabla}_\infty$; let $\wt D_n(\tth)$ stand for the
push-forward of $D_n(\tth)$ under this map. Then $PD(\tth)$ appears as the weak
limit of the measures $\wt D_n(\tth)$ as $n\to\infty$.

\medskip

(c) {\it Projection of a product measure\/} \cite{Ver77a},
\cite[\S4.11]{Arr03}. Consider the infinite-dimensional simplex
$$
\overline{\De}_\infty:=\big\{(u_1,u_2,\dots)\in[0,1]^\infty\;\big|\;
u_1+u_2+\dots\le1\big\}.
$$
The triangular transformation of coordinates $v=(v_1,v_2,\dots)\mapsto
u=(u_1,u_2,\dots)$ given by
$$
u_1=v_1; \qquad u_n=v_n(1-v_1)\dots(1-v_{n-1}),  \quad n\ge2,
$$
maps the cube $[0,1]^\infty$ onto the simplex $\overline{\De}_\infty$. This map
is almost one-to-one: it admits the inversion $u\mapsto v$,
\begin{equation}\label{u-to-v}
v_1=u_1; \qquad v_n=\frac{u_n}{1-u_1-\dots-u_{n-1}}\,, \quad n\ge2,
\end{equation}
which is well defined provided that all the partial sums of the series
$u_1+u_2+\dots$ are strictly less than 1.

Next, the rearrangement of coordinates in weakly decreasing order determines a
projection $\overline{\De}_\infty\to\overline{\nabla}_\infty$.

Denoting by $B(\tth)$ the probability measure on $[0,1]^\infty$ obtained as the
product of infinitely many copies of the measure $\tth(1-t)^{\tth-1}dt$ on
$[0,1]$, the Poisson--Dirichlet distribution $PD(\tth)$ coincides with the
push-forward of $B(\tth)$ under the composition map
$[0,1]^\infty\to\overline{\De}_\infty\to\overline{\nabla}_\infty$.

\medskip

(d) {\it Characterization via correlation functions\/} \cite[\S3]{Wat76}.
Removing possible 0's from a sequence $x\in\overline{\nabla}_\infty$ one may
interpret it as a locally finite point configuration on the semi-open interval
$(0,1]$. This allows one to interpret any probability measure on
$\overline{\nabla}_\infty$ as a random point process on $(0,1]$ (see
\cite{Bor10} for basic definitions). It turns out that the correlation
functions of the point process associated to $PD(\tth)$ have a very simple
form:
$$
\rho_n(u_1,\dots,u_n)=\begin{cases}
\dfrac{\tth^n(1-u_1-\dots-u_n)^{\tth-1}}{u_1\dots
u_n}, & \sum\limits_{i=1}^n u_i< 1;\\
0, & \text{\rm otherwise}. \end{cases}
$$
This provides one more characterization of $PD(\tth)$.

\medskip

The literature devoted to the Poisson--Dirichlet distributions  and their
various connections and applications is very large. The interested reader will
find a rich material in \cite{Arr03}, \cite{Ver72}, \cite{Wat76},
\cite{Ver77a}, \cite{Ver78}, \cite{Ign82}, \cite{Pit97}, \cite{Hol01}
\cite{Kin75}.

Note that $PD(\tth)$ describes the asymptotics of the {\it large\/}
cycle-lengths of random permutations with law $P^\n_\tth$ (namely, the $i$th
coordinate $x_i$ on $\overline\nabla_\infty$ corresponds to the $i$th largest
cycle-length scaled by the factor of $1/n$). The literature also contains
results concerning the asymptotics of other statistics on random permutations,
for instance, small cycle-lengths and the number of cycles \cite{Arr03}.

\section{The Plancherel measure}\label{3}

\subsection{Definition of the Plancherel measure}\label{3-1}

Partitions parameterize not only the conjugacy classes in the symmetric groups
but also their irreducible representations. So far we focused on the conjugacy
classes, but now we will exploit the connection with representations. It is
convenient to identify partitions of $n$ with Young diagrams containing $n$
boxes. The set of such diagrams will be denoted as $\Y_n$. Given a diagram
$\la\in\Y_n$, let $V_\la$ denote the corresponding irreducible representation
of $S_n$ and $\dim\la$ its dimension. In particular, the one-row diagram
$\la=(n)$ and the one-column diagram $\la=(1^n)$ correspond to the only
one-dimensional representations, the trivial and the sign ones. Note that the
symmetry map $\Y_n\to\Y_n$ given by transposition $\la\mapsto\la'$ amounts to
tensoring $V_\la$ with the sign representation, so that $\dim\la'=\dim\la$.
\footnote{In the context of conjugacy classes the operation of transposition
has no natural interpretation.}

By virtue of Burnside's theorem,
\begin{equation}\label{burnside}
\sum_{\la\in\Y_n}(\dim\la)^2=n!.
\end{equation}
This suggests the definition of a probability distribution $M^\n$ on $\Y_n$:
\begin{equation}\label{planch}
M^\n(\la):=\frac{(\dim\la)^2}{n!}\,, \qquad \la\in\Y_n\,.
\end{equation}
Following \cite{Ver77b}, one calls $M^\n$ the {\it Plancherel measure\/} on
$\Y_n$.

In purely combinatorial terms, $\dim\la$ equals the number of {\it standard
tableaux\/} of shape $\la$ \cite[Section 2.5]{Sag01}. Several explicit
expressions for $\dim\la$ are known, see \cite[Sections 3.10, 3.11]{Sag01}.

\subsection{Limit shape and Gaussian fluctuations}\label{3-2}

\begin{figure}[h]
\begin{center}
 \scalebox{0.5}{\includegraphics{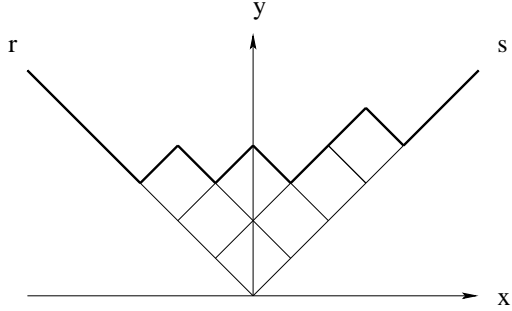}}
\end{center}
\caption{The function $y=\la(x)$ for the Young diagram $\la=(4,2,1)$.
\label{fig}}
\end{figure}

We view each $\la\in\Y_n$ as a plane shape, of area $n$, in the $(r,s)$ plane,
where $r$ is the row coordinate and $s$ the column coordinate. In new
coordinates $x=s-r$, $y=r+s$, the boundary $\partial\la$ of the shape
$\la\subset\R_+^2$ may be viewed as the graph of a continuous piecewise linear
function, which we denote as $y=\la(x)$. Note that $\la'(x)=\pm1$, and $\la(x)$
coincides with $|x|$ for sufficiently large values of $|x|$. The area of the
shape $|x|\le y\le\la(x)$ equals $2n$. (See Figure \ref{fig}.)

Assuming $\la$ to be the random diagram distributed according to the
Plan\-che\-rel measure $M^\n$, we get a random ensemble $\{\la(\,\cdot\,)\}$ of
polygonal lines. We will describe the behavior of this ensemble as
$n\to\infty$.

Informally, the result can be stated as follows: Let $y=\bar\la(x)$ be obtained
from $y=\la(x)$ by shrinking along both the $x$ and $y$ axes with coefficient
 $\sqrt n$,
$$
\bar\la(x)=\frac1{\sqrt n}\,\la(\sqrt n\,x);
$$
then we have
\begin{equation}\label{asympt}
\bar\la(x)\quad\approx\quad \Om(x)+\frac2{\sqrt n}\,\De(x), \qquad n\to\infty,
\end{equation}
where $y=\Om(x)$ is a certain nonrandom curve coinciding with $y=|x|$ outside
$[-2,2]\subset\R$, and $\De(x)$ is a generalized Gaussian process. Let us
explain the exact meaning of \eqref{asympt}.

First of all, the purpose of the scaling $\la(\,\cdot\,)\to \bar\la(\,\cdot\,)$
is to put the random ensembles with varying $n$ on the same scale: note that
the area of the shape $|x|\le y\le\bar\la(x)$ equals 2 for any $n$.

The function $y=\Om(x)$ is given by two different expressions depending on
whether or not $x$ belongs to the interval $[-2,2]\subset\R$:
$$
\Om(x)=\begin{cases} \frac2\pi(x\arcsin\tfrac x2 + \sqrt{4-x^2}), & |x|\le 2\\
|x|, & |x|\ge2
\end{cases}
$$

In the first approximation, the asymptotic relation \eqref{asympt} means
concentration of the random polygonal lines $y=\bar\la(x)$ near a limit curve.
The exact statement (see \cite{Log77}, \cite{Ver77b}, and also \cite{Iva02})
is:

\begin{theorem}[Law of large numbers]\label{omega}
For each $n=1,2,\dots$, let $\la\in\Y_n$ be the random Plancherel diagram and
$\bar\la(\,\cdot\,)$ be the corresponding random curve, as defined above. As
$n\to\infty$, the distance in the uniform metric between $\bar\la(\,\cdot\,)$
and the curve $\Om(x)$ tends to $0$ in probability:
$$
\lim_{n\to\infty}M^\n\big\{\la\in\Y_n\;\big|\;
\sup_{x\in\R}|\bar\la(x)-\Om(x)|\le\varepsilon\big\}=1, \qquad
\forall\varepsilon>0.
$$
\end{theorem}

The second term in the right--hand side of \eqref{asympt} describes the
fluctuations around the limit curve. The Gaussian process $\De(x)$ can be
defined by a random trigonometric series on the interval $[-2,2]\subset\R$, as
follows. Let $\xi_2,\xi_3,\dots$ be independent Gaussian random variables with
mean 0 and variance 1, and set $x=2\cos\theta$, where $0\le\theta\le\pi$. Then
$$
\De(x)=\De(2\cos\theta)
=\frac1\pi\,\sum_{k=2}^\infty\frac{\xi_k}{\sqrt k}\,\sin(k\theta), \qquad
x\in[-2,2].
$$
This is a {\it generalized process\/}, meaning that its trajectories are not
ordinary functions but generalized ones (i.e., distributions). In other words,
it is a Gaussian measure on the space of distributions supported by $[-2,2]$.
For any smooth test function $\varphi$ on $\R$, the smoothed series
$$
\frac1{\pi}\,\sum_{k=2}^\infty \frac{\xi_k}{\sqrt k}\, \int_{-2}^2
\sin(k\theta)\varphi(x)dx, \qquad \theta=\arccos(x/2),
$$
converges and represents a Gaussian random variable. However, the value of
$\De(x)$ at a point $x$ is not defined.

More precisely, the result about the Gaussian fluctuations looks as follows:

\begin{theorem}[Central limit theorem for global fluctuations]\label{global}
Let, as above, $\{\bar\la(\,\cdot\,)\}$ be the random ensemble governed by the
Plancherel measure $M^\n$, and set
$$
\De_n(x)=\frac{\sqrt n}2(\bar\la(x)-\Om(x)), \qquad x\in\R.
$$
For any finite collection of
polynomials $\varphi_1(x),\dots,\varphi_m(x)$, the joint distribution of the
random variables
\begin{equation}\label{int}
\int_\R \varphi_i(x)\De_n(x)dx, \qquad 1\le i\le m
\end{equation}
converges, as $n\to\infty$, to that of the Gaussian random variables
$$
\int_\R\varphi_i(x)\De(x) dx, \qquad 1\le i\le m.
$$
\end{theorem}

This result is due to Kerov \cite{Ker93a}; a detailed exposition is given in
\cite{Iva02}.

Note that for any diagram $\la$, the function $\bar\la(x)-\Om(x)$ vanishes for
$|x|$ large enough, so that the integral in \eqref{int} makes sense.

The theorem implies that the normalized fluctuations $\De_n(x)$, when
appropriately smoothed, are  of finite order. This can be rephrased by saying
that, in the $(r,s)$ coordinates, the global fluctuations of the boundary
$\partial\la$ of the random Plancherel diagram $\la\in\Y_n$ in the direction
parallel to the diagonal $r=s$ have finite order.

A different central limit theorem is stated in \cite{BS}: that result describes
fluctuations at points (so that there is no smoothing); then an additional
scaling of order $\sqrt{\log n}$ along the $y$-axis is required.

Theorem \ref{omega} should be compared to a similar concentration result for
spectra of random matrices (convergence to Wigner's semicircle law). A
similarity between the two pictures becomes especially convincing in view of
the fact (discovered in \cite{Ker93b}) that there is a natural transform
relating the curve $\Om$ to the semicircle law. As for Theorem \ref{global}, it
has a strong resemblance to the central limit theorems for random matrix
ensembles, established in \cite{Dia94}, \cite{Joh98}.

Biane \cite{Bia01} considered a modification of the Plancherel measures $M^\n$
related to the Schur-Weyl duality and found a one-parameter family of limit
curves forming a deformation of $\Om$.

\subsection{The poissonized Plancherel measure as a determinantal
process}\label{3-3}

Let $\Y=\Y_0\cup\Y_1\cup\dots$ be the countable set of all Young diagrams
including the empty diagram $\varnothing$. To each $\la\in\Y$ we assign an
infinite subset $\mathcal L(\la)$ on the lattice $\Z':=\Z+\frac12$ of
half-integers, as follows
\begin{equation}\label{L}
\mathcal L(\la)=\big\{\la_i-i+\tfrac12\;\big|\; i=1,2,\dots\big\}.
\end{equation}
We interpret $\mathcal L(\la)$ as a particle configuration on the nodes of the
lattice $\Z'$ and regard the unoccupied nodes $\Z'\setminus\mathcal L(\la)$ as
{\it holes\/}. In particular, the configuration $\mathcal L(\varnothing)$ is
$\Z'_-:=\{\dots, -\frac32,-\frac12\}$ and the corresponding holes occupy
$\Z'_+:=\{\frac12,\frac32,\dots\}$. In this picture, appending a box to a
diagram $\la$ results in moving a particle from $\mathcal L(\la)$ to the
neighboring position on the right. Thus, growing $\la$ from the empty diagram
$\varnothing$ by consecutively appending a box can be interpreted as a passage
from the configuration $\Z'_-$ to the configuration $\mathcal L(\la)$ by moving
at each step one of the particles to the right by 1.

The configurations $\mathcal L(\la)$ are precisely those configurations for
which the number of particles on $\Z'_+$ is finite and equal to the number of
holes on $\Z'_-$. Note also that transposition $\la\to\la'$ translates as
replacing particles by holes and vice versa, combined with the reflection map
$x\to -x$ on $\Z'$.

The {\it poissonized\/} Plancherel measure with parameter $\nu>0$ \cite{Bai99}
is a probability measure $M_\nu$ on $\Y$, which is obtained by mixing together
the measures $M^\n$ (see \eqref{planch}), $n=0,1,2,\dots$, by means of a
Poisson distribution on the set of indices $n$:
\begin{equation*}
M_\nu (\la) = e^{-\nu}\frac{\nu^{|\la|}}{|\la|!}\, M^{(|\la|)}(\la) = e^{-\nu}
\nu^{|\la|} \, \left(\frac{\dim\la}{|\la|!} \right)^2 \,, \qquad \la\in\Y
\end{equation*}
(see also \cite[Sect. 1.6]{Bor10}).

\begin{theorem}\label{Bessel}
Under the correspondence $\la\to\mathcal L(\la)$ defined by \eqref{L}, the
poissonized Plancherel measure $M_\nu$ turns into the determinantal point
process on the lattice $\Z'$ whose correlation kernel is the discrete Bessel
kernel.
\end{theorem}

About determinantal point processes in general, see \cite{Bor10}. The discrete
Bessel kernel is written down in \cite[Sect. 11.6]{Bor10} (replace there $\tth$
by $\sqrt\nu$). Note that it is a projection kernel. Theorem \ref{Bessel} was
obtained in \cite{Joh01a} and (in a slightly different form) in \cite{Bor00b}.
Johansson's approach \cite{Joh01a} is also discussed in his note \cite{Joh01b}
and the expository paper \cite{Joh05}.

\subsection{The bulk limit}\label{3-4} (See \cite{Bor00b}.)
Fix $a\in(-2,2)$. Recall that the point $(a,\Om(a))$ on the limit curve
$y=\Om(x)$ (see Section \ref{4-2}) corresponds to the intersection of the
boundary $\partial \la$ of the typical large Plancherel diagram $\la\in\Y_n$
with the line $j-i=a\sqrt n$.  The next result describes the asymptotic
behavior of the boundary $\partial\la$ near this point.

\begin{theorem}\label{sine}
Assume that $n\to\infty$ and $x(n)\in\Z'$ varies together with $n$ in such a
way that $x(n)/\sqrt n\to a\in(-2,2)$. Let $\la\in\Y_n$ be the random diagram
with law $M^\n$ given by \eqref{planch} and let $X_n$ be the random particle
configuration on $\Z$ obtained from the configuration $\mathcal L(\la)$ defined
by \eqref{L} under the shift $x\mapsto x-x(n)$ mapping $\Z'$ onto $\Z$. Then
$X_n$ converges to a translation invariant point process on $\Z$, with the
correlation kernel
$$
S^a(k,l)=\begin{cases}\dfrac{\sin(\arccos(a/2)(k-l))}{\pi(k-l)}, & k,l\in\Z,
\quad k\ne l;\\ \dfrac{\arccos(a/2)}\pi, & k=l.\end{cases}
$$
\end{theorem}

The kernel $S^a(k,l)$ is called the {\it discrete sine kernel\/}. It is a
projection kernel and should be viewed as a lattice analogue of the famous sine
kernel on $\R$ originated in random matrix theory. Like the sine kernel, the
discrete sine kernel possesses a universality property \cite{Bai07}.

Theorem \ref{sine} is derived from Theorem \ref{Bessel}: Let $J^\nu(x,y)$
denote the discrete Bessel kernel; one shows that
$$
\lim_{\nu\to\infty}J^\nu(x(\nu)+k,\, x(\nu)+l)=S^a(k,l), \qquad x(\nu)\in\Z',
\quad x(\nu)\sim a\sqrt\nu,
$$
and then one applies a depoissonization argument to check that the large-$n$
limit and the large-$\nu$ limit are equivalent.

\subsection{The edge limit}\label{3-5}

\begin{theorem}\label{airy}
Let $\la=(\la_1,\la_2,\dots)\in\Y_n$ be distributed according to the $n$th
Plancherel measure $M^\n$ given by \eqref{planch}. For any fixed $k=1,2,\dots$,
introduce real-valued random variables $u_1,\dots,u_k$ by setting
\begin{equation}\label{scaling}
\la_i=2n^{1/2}+u_in^{1/6}, \qquad i=1,\dots,k.
\end{equation}
Then, as $n\to\infty$, the joint distribution of $u_1,\dots,u_k$ converges to
that of the first $k$ particles in the Airy point process.
\end{theorem}

Recall (\cite[Sect. 1.9]{Bor10}) that the {\it Airy point process\/} is a
determinantal process on $\R$ living on point configurations $(u_1>u_2>\dots)$
bounded from above; it is determined by the {\it Airy correlation kernel\/},
which is a projection kernel on $\R$.

The Airy point process arises in the edge limit transition from a large class
of random matrix ensembles. It turns out that it also describes the limit
distribution of a few (appropriately scaled) largest rows of the random
Plancherel diagram.

Due to symmetry of $M^\n$ under transposition $\la\to\la'$, the same result
holds for the largest column lengths as well.

Already the simplest case $k=1$ of Theorem \ref{airy} is very interesting,
especially  because of its connection to longest increasing subsequences in
random permutations (see Section \ref{3-6} below). The claim for $k=1$ was
first established by Baik, Deift, and Johansson \cite{Bai99}; then they proved
the claim for $k=2$, \cite{Bai00}. Their work completed a long series of
investigations and at the same time opened the way to generalizations. The
general case of Theorem \ref{airy} is due to Okounkov \cite{Oko00}; note that
his approach is very different from that of \cite{Bai99}, \cite{Bai00}. Shortly
afterwards, the theorem was obtained by yet another method in independent
papers \cite{Bor00b} and \cite{Joh01a}, by using Theorem \ref{Bessel} as an
intermediate step. Note that once one knows Theorem \ref{Bessel}, the precise
form of the scaling \eqref{scaling} can be guessed by a simple argument, see
\cite{Ols08}.

\subsection{Longest increasing subsequences}\label{3-6}

Given a permutation $s\in S_n$, let $L_n(s)$ stand for the length of the {\it
longest increasing subsequence\/} in the permutation word $\wt s:=s(1)s(2)\dots
s(n)$. \footnote{An increasing subsequence in $\wt s$ is a subword $s(i_1)\dots
s(i_k)$ such that $i_1<\dots<i_k$ and $s(i_1)<\dots<s(i_k)$.} Under the uniform
distribution on $S_n$, $L_n$ becomes a random variable. In the sixties, S.~Ulam
raised the question about its asymptotic properties as $n\to\infty$. This
seemingly rather particular problem turned out to be surprisingly deep (about
the history of the problem and many related results, see \cite{Bai99} and the
survey papers \cite{Ald99}, \cite{Dei00}, \cite{Sta07}). The next claim relates
$L_n$ to the Plancherel measure $M^\n$:

\begin{theorem}\label{ulam}
The distribution of $L_n$ under the uniform measure on $S_n$ coincides with the
distribution of $\la_1$, the first row length of the random Young diagram
$\la\in\Y_n$ with law $M^\n$.
\end{theorem}

This result is obtained with the help of the {\it Robinson-Schensted
correspondence\/}, which establishes an explicit bijection $RS:
s\leftrightarrow (\PP,\QQ)$ between permutations $s\in S_n$ and couples
$(\PP,\QQ)$ of standard tableaux of one and the same shape $\la\in\Y_n$. The
bijection $RS$ is described in detail in many textbooks, e.g., \cite{Ful97} and
\cite{Sag01}. The latter book also contains an elegant geometric interpretation
of $RS$ due to Viennot. By the very definitions, the push-forward under $RS$ of
the uniform measure on $S_n$ is $M^\n$. A nontrivial fact is that under this
bijection, $L_n(s)=\la_1$.

By virtue of Theorem \ref{ulam}, the Ulam problem is completely solved by the
$k=1$ case of Theorem \ref{airy} discussed above: the limit distribution of the
scaled random variable $(L_n-2\sqrt n)n^{-1/6}$ is the GUE Tracy-Widom
distribution $F_2$ \cite{Tra94}.

Given a subset $S^*_n\subset S_n$, denote by $L^*_n$ the random variable
$L_n(\,\cdot\,)$ directed by the uniform measure on $S^*_n$. A modification of
the Ulam problem consists in studying the limit distribution of $L^*_n$
(suitably centered and scaled) for subsets $S^*_n$ determined by certain
symmetry conditions imposed on the matrix $[s_{ij}]$ of a permutation $s\in
S_n$. Baik and Rains (see \cite{Bai01} and references therein) showed that in
this way one can get two other Tracy-Widom distributions \cite{Tra96}, $F_1$
and $F_4$, as well as a large family of allied probability distributions
including an interpolation between $F_1$ and $F_4$. These results demonstrate
once again a similarity in asymptotic properties of random permutations and
random matrices. Here is the simplest example from \cite{Bai01}, which shows
that involutions $s=s^{-1}$ in $S_n$ (i.e., symmetric permutation matrices)
model real symmetric matrices:

\begin{theorem}\label{inv}
Take as $S^*_n$ the subset of involutions in $S_n$, and let $L^*_n$ be the
corresponding random variable. Then the limit distribution of $(L^*_n-2\sqrt
n)n^{-1/6}$ is the GOE Tracy-Widom distribution $F_1$.
\end{theorem}

\section{The z-measures and Schur measures}\label{4}

\subsection{The z-measures}\label{4-1}

The identity \eqref{burnside} admits an extension depending on two parameters
$z,z'\in\C$:
\begin{equation*}
\sum_{\la\in\Y_n}(z)_\la(z')_\la(\dim\la)^2=(zz')_n n!,
\end{equation*}
where $(x)_n:=x(x+1)\dots(x+n-1)$ is the Pochhammer symbol and $(x)_\la$ is its
generalization,
\begin{equation*}
(x)_\la:=\prod_{(i,j)\in\la}(x+j-i),
\end{equation*}
the product taken over the boxes $(i,j)$ belonging to $\la$, where $i$ and $j$
stand for the row and column number of a box. The (complex-valued) {\it
z-measure\/} $M^\n_\zz$ on $\Y_n$ assigns weights
$$
M^\n_\zz(\la)=\frac{(z)_\la(z')_\la}{(zz')_n}M^\n(\la)
=\frac{(z)_\la(z')_\la}{(zz')_n}\,\frac{(\dim\la)^2}{n!}
$$
to diagrams $\la\in\Y_n$. This is a deformation of the Plancherel measure
$M^\n$ in the sense that $M^\n_\zz(\la)\to M^\n(\la)$ as $z,z'\to\infty$. In
what follows we assume that the parameters take {\it admissible values\/}
meaning that $(z)_\la(z')_\la\ge0$ for any $\la\in\Y$ and $zz'>0$ (for
instance, one may assume $z'=\bar z\in\C\setminus\{0\}$).  Then $M^n_\zz$ is a
probability measure for every $n$.

The z-measures first emerged in \cite{Ker93c}; they play an important role in
the representation theory of the infinite symmetric group $S_\infty$: Recall
that in Section \ref{2-3} we have mentioned generalized regular representations
$T_z$; it turns out that when $z'=\bar z$, a suitably defined large-$n$ scaled
limit of the z-measures governs the spectral decomposition of $T_z$ into
irreducibles: \cite[\S3]{Bor01a}, \cite{Ols03b}.

The {\it mixed\/} z-measure $M_\zzxi$ on $\Y$ with admissible parameters $(z,
z')$ and an additional parameter $\xi\in(0,1)$ is obtained by mixing up the
z-measures with varying superscript $n$ by means of a negative binomial
distribution on $\Z_+$:
\begin{equation*}
M_\zzxi(\la)=(1-\xi)^{zz'}\frac{(zz')_{|\la|}\xi^{|\la|}}{|\la|!}M^{(|\la|)}_\zz(\la)
=(1-\xi)^{zz'}\xi^{|\la|}(z)_{|\la|}(z')_{|\la|}\left(\frac{\dim\la}{|\la|!}\right)^2,
\end{equation*}
where $\la$ ranges over $\Y$. This procedure is similar to poissonization of
the Plancherel measure and serves the same purpose of facilitating the study of
limit transitions. Note that the poissonized Plancherel measure $M_\nu$ is a
degeneration of $M_\zzxi$ when $z,z'\to\infty$ and $\xi\to0$ in such a way that
$zz'\xi\to\nu$.

\begin{theorem}\label{hypergeom}
Under the correspondence $\la\to\mathcal L(\la)$ defined by \eqref{L}, the
mixed z-measure $M_\zzxi$ turns into a determinantal point process on the
lattice $\Z'$ whose correlation kernel can be explicitly expressed through the
Gauss hypergeometric function.
\end{theorem}

This is a generalization of Theorem \ref{Bessel}. Various proofs have been
given in \cite{Bor00a}, \cite{Bor00c}, \cite{Oko01b}, \cite{Bor06}.

For the lattice determinantal process from Theorem \ref{hypergeom} there are
three interesting limit regimes, as $\xi\to1$, leading to continuous and
discrete determinantal processes:

\smallskip

\begin{itemize}

\item[(1)] Split $\Z'$ into positive and negative parts,
$\Z'=\Z'_+\sqcup\Z'_-$. Given $\la\in\Y$, let $\mathcal L^\circ(\la)\subset\Z'$
be obtained from $\mathcal L(\la)$ by switching from particles to holes on
$\Z'_-$; then $\mathcal L^\circ(\la)$ is finite and contains equally many
particles in $\Z'_+$ and in $\Z'_-$. Note that this {\it particle/hole
involution\/} does not affect the determinantal property. Next, scale the
lattice $\Z'$ making its mesh equal to small parameter $\epsi=1-\xi$. Letting
$\xi\to1$, one gets in this way from $M_\zzxi$ a determinantal process living
on the punctured real line $\R\setminus\{0\}$. The corresponding correlation
kernel is called the {\it Whittaker kernel\/}, because it is expressed through
the classical Whittaker function. This limit process is of great interest for
harmonic analysis on the infinite symmetric groups. For more detail, see
\cite{Bor00a}, \cite{Ols03b}.

\smallskip

\item[(2)] No scaling, we remain on the lattice. The limit determinantal
process is directed by a diffuse measure on the space $\{0,1\}^{\Z'}$ of all
lattice point configurations, and the limit correlation kernel is expressed
through Euler's gamma function, see \cite{Bor05c}.

\smallskip

\item[(3)] An ``intermediate'' limit regime assuming a scaling. It leads to a
stationary limit process whose correlation kernel is expressed through
trigonometric functions and is a deformation of the sine kernel, see
\cite{Bor05c}.

\end{itemize}

\smallskip

These three different regimes describe the asymptotics of the largest, smallest
and intermediate Frobenius coordinates of random Young diagrams, respectively.

\medskip

\begin{remark}
Note a special role of the quantity $\dim\la/|\la|!$ in the expression
for $M_\zzxi$: this is a Vandermonde-like object, which creates a kind of {\it
log-gas pair interaction\/} between particles from the random configuration
$\mathcal L(\la)$ (about log-gas systems, see \cite{For10a}, \cite{For10b}).
The particle/hole involution $\mathcal L(\la)\to\mathcal L^\circ(\la)$ changes
the sign of interaction between particles on the different sides from 0, so
that we get two kinds of particles which are oppositely charged. Note that in
the first regime, the particle/hole involution is necessary for existence of a
limiting point process. The Whittaker kernel is an instance of a correlation
kernel which is symmetric with respect to an {\it indefinite\/} inner product.
\end{remark}

\subsection{Special instances of z-measures}\label{4-2}

$\phantom{aa}$ (a) {\it Meixner and Laguerre ensembles\/}. Assume
$z=N=1,2,\dots$ and $z'=N+b-1$ with $b>0$; these are admissible values. Then
$M_\zzxi$ is supported by the subset $\Y(N)\subset\Y$ of Young diagrams with at
most $N$ nonzero rows. Under the correspondence
$$
\Y(N)\ni\la\,\mapsto\,(l_1,l_2,\dots,l_N)=(\la_1+N-1,\la_2+N-2,\dots,\la_N)\subset\Z_+,
$$
the measure turns into a random-matrix-type object: the $N$-particle Meixner
orthogonal polynomial ensemble with the discrete weight function
$(b)_l\xi^l/l!$, where the argument $l$ ranges over $\Z_+$ (for generalities
about orthogonal ensembles, see \cite{Kon05}). It follows that for general
values of $(z,z')$, the measure $M_\zzxi$ may be viewed as the result of {\it
analytic continuation\/} of the Meixner ensembles with respect to parameters
$N$ and $b$. This observation is exploited in \cite{Bor06}. In a scaling limit
regime as $\xi\to1$, the $N$-particle Meixner ensemble turns into the $N$-point
Laguerre ensemble; the correlation kernel for the latter ensemble is a
degeneration of the Whittaker kernel, see \cite{Bor00a}.

(b) {\it Generalized permutations\/}. Recall that the Plancherel measure $M^\n$
on $\Y_n$ coincides with the push-forward of the uniform measure on $S_n$ under
the projection $S_n\to\Y_n$ afforded by the Robinson--Schensted correspondence
$RS$ (Section\ref{3-6}). Here is a generalization:

Fix natural numbers $N\le N'$ and replace $S_n$ by the finite set $S^\n_{N,N'}$
consisting of all $N\times N'$ matrices with entries in $\Z_+$ such that sum of
all entries equals $n$. Elements of $S^\n_{N,N'}$ are called {\it generalized
permutations\/}. Knuth's generalization of the Robinson--Schensted
correspondence (the $RSK$ correspondence, see, e.g., \cite[Section 4.1]{Ful97})
provides a projection of $S^\n_{N,N'}$ onto $\Y_n(N):=\Y_n\cap\Y(N)$, the set
of Young diagram with $n$ boxes and at most $N$ nonzero rows. It turns out that
the push-forward of the uniform distribution on $S^\n_{N,N'}$ coincides with
the z-measure $M^\n_{N,N'}$, see \cite{Bor01a}.

(c) {\it A variation\/}. In the same way one can get the mixed z-measure
$M_{N,N',\xi}$ if instead of $S^\n_{N,N'}$ one takes $N\times N'$ matrices
whose entries are i.i.d. random variables, the law being the geometric
distribution with parameter $\xi$.

(d) {\it Random words\/}. Denote by $S^\n_{N,\infty}$ the set of words of
length $n$ in the alphabet $[N]:=\{1,\dots,N\}$. Endowing $S^\n_{N,\infty}$
with the uniform measure we get a model of random words. This model may be
viewed as a degeneration of the model of random generalized permutations (item
(b) above) in the limit $N'\to\infty$ (this explains the notation
$S^\n_{N,\infty}$). The $RSK$ correspondence (or rather its simpler version due
to Schensted) provides a projection $S^\n_{N,\infty}\to\Y_n(N)$ taking random
words to random Young diagrams $\la\in\Y_n(N)$ with distribution
$M^\n_{N,\infty}:=\lim_{N'\to\infty}M^\n_{N,N'}$. Asymptotic properties of
random words are studied in \cite{Tra01} and \cite{Joh01a}. The model of random
words can be generalized by allowing non-uniform probability distributions on
the alphabet (see \cite{Its01} and references therein). As explained in
\cite{Its01}, this more general model is connected to the Schur measure
discussed in Section \ref{4-3} below.

(e) {\it The Charlier ensemble and the Plancherel degeneration\/}.
Poissonization of the measure $M^\n_{N,\infty}$ with respect to parameter $n$
leads to the $N$-particle Charlier ensemble \cite[\S9]{Bor01a}. Alternatively,
it can be obtained as a limit case of the mixed z-measures $M_{N,N',\xi}$. The
poissonized Plancherel measure $M_\nu$ appears as the limit of the mixed
z-measures $M_{z,z',\xi}$ when $z,z'\to\infty$ and $\xi\to0$ in such a way that
$zz'\xi\to\nu$. This fact prompted the derivation of the discrete Bessel kernel
(Theorem \ref{Bessel}) in \cite{Bor00b}. Alternatively, $M_\nu$ can be obtained
through a limit transition from the Charlier or Meixner ensembles; this leads
to another derivation of the discrete Bessel kernel: \cite{Joh01a},
\cite{Joh01b}.

\subsection{The Schur measures}\label{4-3}

Let $\La$ denote the graded algebra of symmetric functions. The Schur functions
$s_\la$, indexed by arbitrary partitions $\la\in\Y$, form a distinguished
homogeneous basis in $\La$. As a graded algebra, $\La$ is isomorphic to the
algebra of polynomials in countably many generators; as these generators, one
can take, for instance, the complete homogeneous symmetric functions
$h_1,h_2,\dots$ where $\deg h_k=k$. One has $s_\la=\det[h_{\la_i-i+j}]$ with
the understanding that $h_0=1$ and $h_k=0$ for $k<0$ (the Jacobi-Trudi
formula); here the order of the determinant can be chosen arbitrarily provided
it is large enough. For more detail, see, e.g., \cite{Sag01}.

Given two multiplicative functionals $\varphi,\psi\colon\La\to\C$, the
corresponding (com\-plex-valued) {\it Schur measure\/} $M_{\varphi,\psi}$ on
$\Y$ is defined by
\begin{equation*}
M_{\varphi,\psi}(\la)=\const^{-1}\varphi(s_\la)\psi(s_{\la}), \quad \la\in\Y,
\qquad \const=\sum_{\la\in\Y}\varphi(s_\la)\psi(s_{\la}),
\end{equation*}
provided that the sum is absolutely convergent (which is a necessary condition
on $\varphi,\psi$). This notion, due to Okounkov \cite{Oko01a}, provides a
broad generalization of the mixed z-measures. Since a multiplicative functional
is uniquely determined by its values on the generators $h_k$, the Schur measure
has a doubly-infinite collection of parameters $\{\varphi(h_k), \psi(h_k);
k=1,2,\dots\}$. In this picture, the z-measures correspond to a very special
collection of parameters
$$
\varphi(h_k)=\xi^{k/2}(z)_k/k!, \quad \psi(h_k)=\xi^{k/2}(z')_k/k!, \qquad
k=1,2,\dots,
$$
and the poissonized Plancherel measure $M_\nu$ appears when
$\varphi(h_k)=\psi(h_k)=\nu^{k/2}/k!$.

As shown in \cite{Oko01a},  Theorem \ref{hypergeom} extends to Schur measures:
if the parameters are such that the measure $M_{\varphi,\psi}$ is nonnegative
(and hence is a probability measure), then it gives rise to a lattice
determinantal point process. Moreover, for the corresponding correlation kernel
one can write down an explicit contour integral representation \cite{Bor00c}.
Such a representation is well suited for asymptotic analysis.

If $\varphi$ and $\psi$ are evaluations of symmetric functions at finitely many
positive variables, the first row $\lambda_1$ can be interpreted as the last
passage percolation time in a suitable directed percolation model on the plane,
see \cite{Joh05}.

\subsection{Some generalizations}\label{4-4}

Kerov \cite{Ker00} generalized the construction of the z-measures
$M^{(n)}_{z,z'}$ by introducing an additional parameter related to Jack
polynomials. This new parameter is similar to the $\beta$ parameter in random
matrix ensembles \cite{For10b}. In particular, the Plancherel measure
$M^{(n)}$, which is a limit case of the z-measures, also allows a
$\beta$-deformation \cite{Ker00}, \cite{Oko05}, \cite{Oko06}. The ordinary
z-measures correspond to the special value $\be=2$, and in the limit $\be\to0$
the beta z-measures degenerate to the measures \eqref{ESF} derived from the
Ewens measures, see \cite[Section 1.2]{Ols10}. Thus, the $\be$ parameter
interpolates between the models of Section \ref{2} and those of Sections
\ref{3}-\ref{4}, as has been pointed out in the end of Section \ref1. Note also
that replacing the Schur functions by the Jack symmetric functions leads to a
natural $\be$-deformation of the Schur measures.

As in random matrix theory, the value $\be=2$ is a distinguished one, while in
the general case $\be>0$ the situation is much more complex. Some results for
$\be\ne2$ can be found in \cite{Bor05b}, \cite{Ful04}, \cite{Ols10},
\cite{Str10a}, \cite{Str10b}.

In a somewhat different direction, one can define natural analogues of the
Plancherel measure and Schur measures for shifted Young diagrams (equivalently,
strict partitions): \cite{Tra04}, \cite{Mat05}. This theory is related to
Schur's Q-functions (a special case of Hall--Littlewood symmetric functions
that appears in the theory of projective representations of the symmetric
group). Surprisingly enough, a natural analogue of the z-measures for shifted
diagrams, discovered by Borodin and recently studied in \cite{Pet10}, seems to
be not related to Schur's Q-functions.

Finally, note that there are many points of contact between the results
described in this chapter and Fulman's work on ``random matrix theory over
finite fields'', see his survey \cite{Ful01} and references therein.

\newpage

{\sc Acknowledgements}: This work was supported by the RFBR grant 08-01-00110.
I am grateful to Jinho Baik, Alexei Borodin, Alexander Gnedin, and the referee
for helpful comments.

\end{document}